\documentclass{amsart}

\usepackage{amssymb}
\usepackage{hyperref}
\usepackage[symbol]{footmisc}

\begin{document}

\newtheorem{theorem}[equation]{Theorem}
\newtheorem{proposition}[equation]{Proposition}
\newtheorem{lemma}[equation]{Lemma}
\theoremstyle{definition}
\newtheorem{definition}[equation]{Definition}
\newtheorem*{notation*}{Notation}
\newtheorem{remark}[equation]{Remark}

\numberwithin{equation}{section}

\title[The tame symbol for curves]{Bounding generators for the kernel and cokernel of the tame symbol for curves}

\author{Rob de Jeu}\address{Faculteit der B\`etawetenschappen\\Afdeling Wiskunde\\Vrije Universiteit Amsterdam\\De Boelelaan 1111\\1081 HV Amsterdam\\The Netherlands}

\begin{abstract}
Let $ C $ be a regular, irreducible curve that is projective over a field.
We obtain bounds in terms of the arithmetic genus of~$ C $ for the generators
that are required for the cokernel of the tame symbol, as well
as, under a simplifying assumption, its kernel.
We briefly discuss a potential application to Chow groups.
\end{abstract}

\dedicatory{Dedicated to the memory of Jaap Murre, with admiration and gratitude.}



\maketitle

\def\d{\delta}
\def\e{\varepsilon}

\def\E{{\mathcal E}}
\def\O{{\mathcal O}}

\def\P{\mathbb P}
\def\Q{\mathbb Q}
\def\Z{{\mathbb Z}}

\def\tf{\tilde f}
\def\tS{\widetilde S}

\def\rightiso{\xrightarrow{\cong}}

\def\cl{^{(1)}}

\def\ord{\textup{ord}}
\def\pol #1 {k[x]_{\le #1}}

\section{Introduction}

Let $ F $ be a field.
By a theorem of Matsumoto (see \cite[Theorem~11.1]{mil71}) there is a presentation of $ K_2(F) $
as~$ F^* \otimes_\Z F^* / \langle a \otimes (1-a), a \in F, a \neq 0,1 \rangle $,
where $ \langle\cdots \rangle $ stands for the subgroup of $ F^* \otimes_\Z F^* $ generated by the indicated elements.
The class of $ a \otimes b $ is denoted $ \{ a, b \} $ (which
is called a symbol), so that $ K_2(F) $ is an Abelian group, written additively, with
generators $ \{ a , b \} $ for $ a $ and $ b $ in $ F^* $, and relations  
\begin{alignat*}{1}
& \{ a_1a_2 ,  b \} = \{ a_1 ,  b \} + \{ a_2 ,  b \} 
\\
& \{ a , b_1b_2 \} = \{ a , b_1 \} + \{ a , b_2 \} 
\\
& \{ a , 1-a \} = 0 \text{ if } a \text{ is in } F , \ a \neq 0, 1 .
\end{alignat*}
These relations imply that
$ \{a,b\} + \{b,a\} = \{ c , -c \} = 0 $ for $ a $,  $ b $ and $ c $ in $ F^* $.          

If $ a+b+c=0 $ with $ a $ and $ b $ in $ F^* $
and $ c $ in $ F $, then $ \{a,b\} = 0 $ in $ K_2(F) $ if $ c=0 $, and 
$ \{-a/c,-b/c\} = 0 $ (or equivalently, 
$ \{a,-b\} + \{b,-c\} + \{c,-a\} =  \{-1,-1\} $) if $ c \ne 0 $.
Such identities arise in division with remainder, which can be used
to determine~$ K_2(\Q) $ (see the proof of \cite[Lemma~11.7]{mil71}).

One can similarly use division with remainder in $ k[x] $, where $ k $ is a field
and $ x $ a variable, to compute the kernel and cokernel of the tame symbol
for $ F = k(x) $.  One way of obtaining division with remainder here is by observing that, with $ \pol n $ the~$ k $-vector
space of polynomials in $ k[x] $ of degree at most $ n $, for $ f(x) \ne 0 $ of degree~$ d $ and $ g(x) \ne 0  $
of degree $ d' \le d $, the kernel of the $ k $-linear map
\begin{alignat*}{1}
\pol 0 \times \pol {d-d'} \times \pol {d'-1} & \to \pol d 
\\
(c,q(x), r(x)) & \mapsto c f(x) - q(x) g(x) -r(x)
\end{alignat*}
is non-trivial because of dimensions.  For $ (c,q(x),r(x)) \ne (0,0,0) $ in the kernel
one must have $ c \ne 0 $ because of degrees, so that the kernel has dimension~1.

This argument uses dimensions of Riemann-Roch spaces,
and similar arguments can be given for a 
regular, irreducible curve~$ C $ that is projective over a field~$ k $,
with function field $ F = k(C) $.
Used appropriately, it leads to relations of the form~$ a+b+c=0 $ with $ a $, $ b $ in $ F^* $
and $ c $ in $ F $ that allow rewriting in $ K_2(F) $.
From this, we obtain bounds on the generators that are needed
for the kernel and cokernel of the tame symbol
\begin{equation} \label{TAME}
T : K_2(F) \to \coprod_{ P \in C\cl } k(P)^*
\,.
\end{equation}
Here~$ C\cl $ is the set of codimension~1 (i.e., closed) points of~$ C $,
and the component~$ T_P $ for such a point~$ P $ is given by~$ T_P( \{f, g\} ) = (-1)^{\ord_P(f) \ord_P(g)} \frac{f^{\ord_P(g)}}{g^{\ord_P(f)}} |_P $.
In order to simplify notation, we denote~$ C\cl $ by~$ C $ from
now on.
\def\cl{}
\medskip

The proofs of Proposition~\ref{introcoker} as well as Theorems~\ref{introker1}
and~\ref{introker2} below are essentially algorithmic, so in principle they can be
applied as such to explicit elements.
We also emphasise that the results hold for any field~$ k $
(cf.~\cite[\S5]{Ba-Ta}). In fact, for~$ k $ algebraically closed 
they have limited content because then
in Proposition~\ref{introcoker} we have~$ \tS = C $ unless~$ t = 0 $
and~$ e = 1 $, for which~$ \tS = |D| $,
Theorem~\ref{introker1} is vacuous unless~$  t = 0 $ and~$ d = 1 $,
and~$ S = C $ in Theorem~\ref{introker2}.
Instead, the paper was motivated by the following.

Let~$ E $ be an elliptic curve over a number field~$ k $,
and~$ \E $ a flat, regular, proper model of it over the ring of
algebraic integers of~$ k $.
Then the image of~$ K_2(\E) $ in~$ K_2(k(E)) $ under localisation is contained in the kernel
of the tame symbol, is independent of~$ \E $ (see \cite[Proposition~4.1]{Li-dJ}
and its proof),
and is expected to be finitely generated as a consequence of
a conjecture by Bass (see \cite[Conjecture~IV.6.8]{Kbook}).
We do not approach this last statement, but Theorems~\ref{introker1}
and~\ref{introker2}
show that the kernel of the tame symbol is
contained in a more manageable subgroup of~$ K_2(k(E)) $.

Suppose that~$ C $ has arithmetic genus%
\footnote[2]{We choose this unusual notation so that we can use~$ f $, $ g $ and~$ h $ for functions.}
$ t $ and that~$ H^0(C, \O_C) = k $. It follows that~$ t = \dim_k(H^1(C, \O_C)) \ge 0 $.
One easily checks that~$ \deg(P) \le ne $ for~$ P $ in~$ \tS \setminus |D| $
in the proposition below, so that it provides an upper
bound on the generators that suffice for the cokernel of
the tame symbol.

\begin{proposition} \label{introcoker}
Let~$ C $ and~$ k $ be as above.
For a divisor~$ D $ on~$ C $ of degree~$ e \ge 1 $, 
let~$ n $ be the smallest integer such that~$ n e > 2 t +  e - 2 $,
and let
\begin{equation*}
\tS = |D| \cup \{ P \text { in } C \text{ with } H^0(C , \O_C(nD-P)) \ne \{0\} \}
\,.
\end{equation*}
Then the map
\begin{equation*}
K_2(F) \to \coprod_{P \notin \tS} k(P)^*
\end{equation*}
induced by the tame symbol is surjective, and
the cokernel of the tame symbol
is generated by the image of $ \coprod_{P \in \tS} k(P)^* $.
\end{proposition}

Its proof, which we give in Section~\ref{cokersection}, is not
that complicated in this generality. But although one can undoubtedly
prove results on the kernel of the tame symbol similar to the next two theorems in the
same generality, due to the many technical lemmas involved, this would
probably not provide an optimal result. In order to not obscure
the structure of the proofs, we therefore assume that~$ C $ has a~$ k $-rational
point~$ O $
(which implies~$ H^0(C, \O_C) = k $) and effectively use~$ D = (O) $.

We fix more notation.
For an integer $ d $, let $ C_d\cl = \{ P \text{ in } C\cl \text{ with } \deg(P) = d  \} $,
and put~$ C_{\le d}\cl = \cup_{e \le d} \, C_e\cl $
as well as~$ C_{\ge d}\cl = \cup_{e \ge d} \, C_e\cl $.
Also, let~$ F_d $ the subgroup of~$ K_2(F) $ generated by all symbols $ \{f, g\} $ with~$ |(f)| \cup |(g)| \subseteq C_{\le d}\cl $.

\medskip

As a first step, in Section~\ref{degreesection} we prove the following.

\begin{theorem} \label{introker1}
Let~$ C $, $ k $ and~$ O $ be as above.
The tame symbol induces isomorphisms
\begin{equation*}
F_d / F_{d-1} \rightiso \coprod_{O \ne P \in C_d\cl} k(P)^*
\quad
\text{ and }
\quad
K_2(F) / F_{d-1} \rightiso \coprod_{O \ne P \in C_{\ge d}\cl} k(P)^*
\end{equation*}
for $ d \ge 3t+1 $.
In particular, the tame symbol and its restriction
\begin{equation*}
  F_{d-1}   \to  \coprod_{P \in \{O\} \cup C_{\le d-1}\cl} k(P)^*
\end{equation*}
yield the same kernel and cokernel.
\end{theorem}

Theorem~\ref{introker1} gives a full (and classical) description
of~$ K_2(F) $ for~$ t=0 $ because we can take~$ d = 1 $.
In particular, the kernel of the tame symbol is~$ K_2(k) $
(which one easily checks injects into~$ K_2(F) $ here),
and the cokernel is isomorphic to~$ k^* $.

For $ t \ge 1 $ we refine it in Section~\ref{RRsection} as follows.
Let~$ L_d \subseteq F_d $ be the subgroup of~$ K_2(F) $ generated by symbols $ \{l,m\} $ for non-zero~$ l,m $ in~$ H^0(C, \O_C(d(O))) $,
and let~$ S = \{ P \text{ in }  C\cl \text{ with }  H^0(C, \O_C(3t(O) - (P))) \ne \{0\} \}  $.

\medskip

\begin{theorem} \label{introker2}
Let~$ C $ and~$ O $ be as above.
For~$ t \ge 1 $, the tame symbol induces an isomorphism
\begin{equation*}
K_2(F) / L_{3t} \rightiso \coprod_{P \notin S} k(P)^*
\,.
\end{equation*}
In particular, the tame symbol and its restriction
\begin{equation*}
  L_{3t}   \buildrel{T}\over{\to}  \coprod_{P \in S} k(P)^*
\end{equation*}
yield the same kernel and cokernel.
\end{theorem}

In this case there are also results more fully like those given
in Theorem~\ref{introker1}, but as they are somewhat more
technical to state we refer the reader to Proposition~\ref{RRprop}.
It is also worth noting that this theorem is easier to prove
for~$ t = 1 $, and we do this as an illustration in Section~\ref{RRsection}
before treating the general case.

\medskip

As an example, suppose $ C $ is an elliptic curve over $ k $ with rational point $ O $.
Let~$ \{1,x\} $ and $ \{1,x,y\} $ be bases of  $ H^0(C, \O_C(2(O))) $ and $ H^0(C, \O_C(3(O))) $ respectively, so that $ x $ and $ y $ satisfy a Weierstrass equation.
Then~$ \tS $ in Proposition~\ref{introcoker} consists of~$ O $
as well as all~$ P \ne O $ with~$ x(P) $ in~$ k $.
Theorem~\ref{introker2} states that the kernel of the tame symbol is contained in
the subgroup
\begin{equation*}
L_3 = \langle \{a_1+b_1x+c_1y , a_2+b_2x+c_2y\} \text{ with  all } a_i, b_i \text{ and } c_i \text{ in } k \rangle
\,.
\end{equation*}
Note that each $ a_i+b_ix+c_iy  $ is a constant or defines a
projective line in~$ \P_k^2 $ if we embed the curve
using  $ H^0(C, \O_C(3(O))) $, motivating the notation~`$ L $'.
Moreover, Theorem~\ref{introker2} is sharp in some sense in this case; see Remark~\ref{sharpremark}.

\medskip

As mentioned above, the author wrote the current paper in order to find smaller subgroups
of~$ K_2(F) $ that still contain the kernel of the tame symbol.
But when he discussed the (then preliminary) results with him, it was Jaap Murre who pointed out that it also has potential applications
to Chow groups.
For this, let~$ X $ be a regular, irreducible variety over a field.
In the Gersten-Quillen spectral sequence for~$ X $,
the Chow group~$ CH^p(X) $ of~$ X $ for codimension~$ p \ge 2 $ occurs as the cokernel of the map~$ \coprod_{x \in X^{(p-1)}} k(x)^* \to \coprod_{x \in X^{(p)}} \Z $.
Here~$ X^{(q)} $ denotes the set of points of~$ X $ of codimension~$ q $,
and the map is essentially the one induced by the divisor maps \cite[Proposition~5.26]{Sri96}.
It vanishes on the image of the incoming map~$ \coprod_{x \in X^{(p-2)}} K_2(k(X)) \to \coprod_{x \in X^{(p-1)}} k(x)^* $
induced by the tame symbols, so that Proposition~\ref{introcoker}
has some potential application
to limiting the relations needed among the generators of the Chow group.
As a simple example, suppose~$ X $ is a surface and admits a surjective morphism $ X \to Y $
to a curve~$ Y $, with as generic fibre a curve~$ C $ over the function field~$ k $
of~$ Y $.
Then $ CH^2(X) $ is the cokernel of~$ \coprod_{x \in X^{(1)}} k(x)^* \to \coprod_{x \in X^{(2)}} \Z $.
Writing~$ X^{(1)} $ as the disjoint union of~$X_h^{(1)}$, consisting of the
points mapping to the generic point of~$ Y $, and its complement~$ X_v^{(1)} $,
by Proposition~\ref{introcoker} one can replace~$ X_h^{(1)} $
by a smaller subset and still obtain~$ CH^2(X) $ as the cokernel.

\medskip

We conclude this section with two more points. The first is
about notation.

\begin{notation*}
(1)
For any subset~$ M $ of~$ C\cl $, we write $ M_d $ for $ M \cap C_d\cl $, and
similarly for~$ M_{\le d} $ and~$ M_{\ge d} $.

(2)
For a vector space $ V $ we let $ V^* = V\!\setminus\!\{0\} $.
\end{notation*}

The second consists of the following observation, which, although
in sometimes hidden ways, is the basis for many of our lemmas.

\begin{lemma} \label{VW}
If $ K/k $ is a finite extension of fields, and $ V $, $ W $ are $ k $-subspaces of $ K $
with~$ \dim_k(V) + \dim_k(W) > [K:k] $, then
$ K^* = \{ v w^{-1} \text{ with $ v $ in $ V^* $ and $ w $ in $ W^* $} \} $.
\end{lemma}

It follows immediately from considering,
for~$ \beta $ in $ K^* $, the dimensions of the~$ k $-subspaces~$ V $ and~$ \beta \cdot W $
of~$ K $.
Alternatively, the $ k $-linear map $ V \times W \to K $ given
by mapping $ (v,w) $ to $ v - \beta w $
must have a non-trivial kernel, which is a formulation closer to 
what is
used in the proofs of most lemmas in Sections~\ref{degreesection}
and~\ref{RRsection}.

\section{The cokernel of the tame symbol} \label{cokersection}

We now prove Proposition~\ref{introcoker}.
Because our assumptions are weaker here than in Sections~\ref{degreesection}
and~\ref{RRsection}, we use slightly different notation.

For~$ P $ in~$ C \setminus |D| $ of degree~$ d \ge t + e $,
the Riemann-Roch space~$ L(\lceil \frac{t-d}e \rceil D + (P)) $ has dimension at least~$ 1 $
(see Example~6.3.18 and Remark 7.3.33 of~\cite{Liu}).
Fixing~$ g \ne 0 $ in it, we obtain~$ D_P = (g) + (P) + \lceil \frac{t-d}e \rceil D \ge 0 $.
Then $ \ord_P(g) \ge -1 $ because~$ P $ is not in~$ |D| $, and
$ \ord_P(g) \ge 0 $ is impossible because otherwise $ (g) + \lceil \frac{t-d}e \rceil D \ge 0 $,
which is impossible due to degrees as~$ \lceil \frac{t-d}e \rceil \le -1 $.
From~$ \lceil \frac{t-d}e \rceil - 1 < \frac{t-d}{e} \le \lceil \frac{t-d}e \rceil $
we see that~$ t \le \deg(D_P) < t + e $.
We then set~$ \tf_P = 1/g $, so that~$ (\tf_P) = (P)  - D_P  + \lceil \frac{t-d}e \rceil D $,
and~$ |(\tf_P)| $ is contained in~$ \{P\} \cup |D| \cup C_{\le t+e-1}\cl $.

For $ P $ in $ C \setminus |D| $ with~$ d = \deg(P) \ge 2 t + 2 e - 1 $,
consider a symbol~$ \{ \tf_P, l \} $ with~$ l $ in $ L(mD)^* $ 
for~$ m e < d $.
Then at a point~$ Q \ne P $ in~$ C \setminus |D| $
with~$ \deg(Q) \ge d $ there will be no contribution under the
tame symbol because~$ Q $ is not in~$ |(\tf_P)| $,
the poles of $ l $ are in $ |D| $,
and a zero of~$ l $ at~$ Q $ would imply~$ (l) + m D \ge (Q) $,
giving the contradiction~$ \deg(Q) \le m e $.
This also shows that~$ L(mD) $ injects into~$ k(P) $ under evaluation at~$ P $,
so that for a symbol~$ \{ \tf_P, l_1 / l_2 \} $ with~$ l_1 $ and~$ l_2 $
in~$ L(md)^* $, by Lemma~\ref{VW} we can get a prescribed contribution at~$ P $
when~$ 2 (m e + 1 -t) > d $, i.e.,~$ 2 m e > d + 2 t - 2 $,
and other contributions only at points of degree less than~$ d $
or in~$ |D| $.
The two conditions on~$ m $ can be
satisfied when~$ 2 d > d + 2 t - 2 + 2 e $, which is
our assumption on~$ d $.

From Riemann-Roch one obtains~$ |D| \cup C_{\le t+e-1}\cl \subseteq \tS \subseteq |D| \cup C_{\le 2 t + 2 e - 2}\cl $,
so using the above symbols for decreasing~$ d \ge 2 t + 2 e - 2 $, we are left
to consider points~$ P $ in~$ C \setminus \tS $ of degree~$ d = t + e, \dots, 2 t + 2 e - 2 $.
For any~$ P $ with~$ \deg(P) \ge t+e $
we have~$ |(\tf_P)| \subseteq \{P\} \cup \tS $, and for~$ l $ in~$ L(nD)^* $ we have~$ |(l)| \subseteq \tS $
by the definition of~$ \tS $.
Because~$ L(nD) $ has dimension at least~$ n e + 1 - t \ge t + e $,
and injects into~$ k(P) $ under evaluation at~$ P $ in~$ C \setminus \tS $,
we see from Lemma~\ref{VW} that for~$ P $ in~$ C \setminus \tS $
of degree~$ t + e, \dots, 2 t + 2 e - 2 $ there are~$ l_1 $ and~$ l_2 $
in~$ L(nD)^* $ such that the tame symbol of~$ \{ \tf_P, l_1 / l_2 \} $
has prescribed image at~$ P $,  and all other contributions in~$ \tS $.

\begin{remark}
One sees similarly that also for each~$  n \ge 3 $
the map
\begin{equation*}
 K_n^M(F) \to \coprod_{P \notin \tS} K_{n-1}^M(k(P)) 
\end{equation*}
induced by the tame symbol on Milnor $ K $-theory is surjective.
\end{remark}

\section{Filtration using degrees} \label{degreesection}

Our goal in this section is to prove Theorem~\ref{introker1}.
We recall that~$ C $ is a regular, irreducible
curve that is projective over a field~$ k $, with function field~$ F=k(C) $,
which has arithmetic genus~$ t \ge 0 $ and on which we fixed
a~$ k $-rational point~$ O $.
We first choose a shorter notation for what will be used many times
below.

\begin{notation*}
For a divisor  $ D $ on $ C $ and an integer $ d $, we let~$ RR_d(D) $ be
the Riemann-Roch space~$ L(d(O)-D) $.  When $ D=0 $ we simply write $ RR_d $.
\end{notation*}
\noindent
Then~$ \dim_k(RR_d(D)) \ge d-\deg(D) + 1-t  $, and equality holds
if~$ d-\deg(D) > 2t-2 $.

We shall define some functions~$ f_P $ and~$ f_P' $ that
allow us to express every element of~$ F^* $ in a controlled
way.%
\footnote[3]{%
If we let~$ C' = C \setminus \{O\} \cup C_{\le t} $,
then the~$ f_P $ are prime elements of the
unique factorisation domain~$ H^0(C' , \O_{C'} ) $.
The same holds in Section~\ref{cokersection} for the~$ \tf_P $
and~$ C' = C \setminus | D| \cup C_{\le t +e - 1} $.}
They play a crucial role in our proofs.
For~$ P \ne O $ in~$ C_{\ge t+1}\cl $ with $ \deg(P) = d $
we can find $ f_P $ in $ F^* $ with 
\begin{equation} \label{fPdef}
 (f_P) = (P) - D_P - (d-t)(O)
\end{equation}
for some~$ D_P \ge 0 $ of degree~$ t $, so that $ f_P $ is in
$ RR_{d-t}(-D_P)^* $.  Namely, $ RR_{t-d}(-(P)) $ has positive dimension, so it contains
$ g \ne 0 $ with $ D_P = (g) + (P) + (t-d)(O) \ge 0 $.  So~$ \ord_P(g) \ge -1 $, and
$ \ord_P(g) \ge 0 $ is impossible as then~$ (g) + (t-d)(O) \ge 0 $, which is
impossible because of degrees. We can therefore take $ f_P = 1/g $.
Fixing $ f_P $ as above for every $ P \ne O $ in $ C_{\ge t+1}\cl $,
we can also find (and fix) a non-zero function $ f_P' $ in~$ RR_{2t}(D_P) $
as this space has positive dimension. Then~$ (f_P') = D_P + D_P' - 2t (O) $ for some $ D_P' \ge 0 $,
also of degree $ t $.

\medskip

The following definition corresponds to that of~$ \tS $ in Proposition~\ref{introcoker}
for~$ D = (O) $.

\begin{definition} \label{Sprimedef}
Let~$ S' = \{ O \} \cup \{ P \text{ in }  C\cl \text{ with }  RR_{2t}((P)) \ne \{0\} \} $.
\end{definition}

Then~$ \{ O \} \cup C_{\le t}\cl \subseteq S' \subseteq \{ O \} \cup  C_{\leq 2t}\cl $,
with the first inclusion a consequence of Riemann-Roch.
In particular, for~$ P $ in $  C_{\ge t+1}\cl \setminus \{O\} $, the only point of~$ |(f_P)| $
that does not lie in~$ S' $ can be~$ P $, and it will be outside of $ S' $ if~$ \deg(P) \ge 2t+1 $.

With those choices for~$ f_P $ we can write any $ f $ in $ F^* $ as
\begin{equation} \label{fPprod}
 f = f' \prod_{O \ne P \in |(f)|_{\ge t+1}} f_P^{\ord_P(f)}
\end{equation}
for a unique~$ f' $, which has~$ |(f')| \subseteq C_{\le t}\cl $.
(If~$ t = 0 $ then~$ f' $ is in $ k^* $.)

\medskip

We recall the subgroups~$ F_d $ and~$ L_d $ from the introduction
in the current notation.

\begin{definition}
For any integer $ d $, we let $ F_d $ be the subgroup of $ K_2(F) $ generated by all symbols $ \{f, g\} $
with $ |(f)| \cup |(g)| $ contained in $ C_{\le d}\cl $, and~$ L_d \subseteq F_d $
the subgroup generated by all symbols~$ \{l, m\} $ with $ l $ and $ m $ in $ RR_d^*$.
\end{definition}

We need a lemma on~$ F^* $ in order to rewrite generators of
such subgroups.

\begin{lemma} \label{3tlemma}
(1)
Suppose $ f $ in $ F^* $ satisfies $ |(f)| \subseteq C_{\le t}\cl $.  Then $ f $
is in $ \langle RR_{3t}^* \rangle $.

(2)
If $ P\ne O $ and  $ d =\deg(P) \ge t+1 $ then $ f_P $ is in $ \langle RR_{d+t}^* \rangle $.

(3)
Let $ P \ne O $ have degree $ d \ge t+1 $.  If $ 2t \le d' \le d+t $
and $ RR_{d'}((P)) \ne \{0\} $ then $ f_P $ is in $ \langle RR_{d'}^* \rangle $.
\end{lemma}

\begin{proof}
(1)
This is clear if~$ t = 0 $, so let~$ t \ge 1 $.
If $ \deg(f) \le t $ then $ RR_{2t}((f)_\infty) \subseteq RR_{2t} $ has dimension at least $ 2t-\deg(f) + 1-t = 1+t-\deg(f) \ge 1 $.
For $ g \ne 0 $ in this space we have $ (g) + (f)  + 2t (O) \ge (f)_0 \ge 0 $, so
$ g f $ is in $ RR_{2t}^* $ and the statement is clear.
If $ \deg(f) > t $ then we use induction on $ \deg(f) $.   Since $ |(f)| \subseteq C_{\le t}\cl $
we can write~$ (f)_* = D_* + E_* $ for $ * = 0, \infty $ using effective divisors $ D_*, E_* $
with support in~$ C_{\le t}\cl $ and~$ t < \deg(D_*) \le 2t $.
Computing dimensions, one sees there exist non-zero $ h_* $ in~$ RR_{t+\deg(D_*)}(D_*) \subseteq RR_{3t} $,
so $ (h_*) = D_* + E_*' - (t+\deg(D_*)) (O) $ for effective~$ E_*' $ of degree $ t $.  Then~$ f h_\infty/h_0 $ has
divisor~$ E_0 + E_\infty' - E_\infty - E_0' + (\deg(D_0)-\deg(D_\infty)) (O) $, and this is supported
in~$ C_{\le t}\cl $. With~$ \deg(E_\infty') < \deg(D_0) $
and~$ \deg(E_0') < \deg(D_\infty) $, both~$ E_0 + E_\infty'$ and~$ E_\infty + E_0' $ are effective
divisors of degree less than~$ \deg(f) $. It follows that the
degree of~$ f h_\infty/h_0 $ is less than that of~$ f $,
and by induction is in~$ \langle RR_{3t}^* \rangle $.

(2)
We note that $ (f_P f_P') = (P) + D_P' - (d+t) (O) $ with~$ D_P' \ge 0 $,
and that $ f_P' $ is in $ RR_{2t}(D_P)^* \subseteq RR_{d+t}^* $.

(3)
If $ g $ in $ RR_{d'}^* $ satisfies~$ g(P) = 0 $, then $ (g) = (P) + D_P'' - d' (O) $
with~$ D_P'' \ge 0 $ and of degree $ d'-d \le t $. Therefore~$ (g/f_P) = D_P + D_P'' - (t+d'-d) (O) $,
hence~$ g/f_P $ is in $ RR_{2t}^* \subseteq RR_{d'}^* $.
\end{proof}

Note that~\eqref{fPprod} and parts (1) and (2) of Lemma~\ref{3tlemma}
together imply~$ F_d \subseteq L_{d+t} $ for~$ d \ge 2t $.
Similarly, they imply that $ F^* $ is generated by $ RR_{3t}^* $
together with all~$ f_P $ for $ \deg(P) \ge 2t+1 $, 
because~$ f_P $ lies in $ \langle RR_{3t}^* \rangle $ if $ \deg(P) = t+1,\dots,2t $.
Hence $ F_d $ for $ d \ge 3t $ is generated by the following
symbols:
\begin{equation*}
\begin{aligned}
{\scriptstyle\bullet}\quad & \{ f_P, f_Q \} \text{ with } 2t+1 \le \deg(Q) \le \deg(P) \le d ,
\\
{\scriptstyle\bullet}\quad & \{ f_P, l \} \text{ with } 2t+1 \le \deg(P) \le d \text{ and } l \text{ in } RR_{3t}^* ,
\text{ and }
\\
{\scriptstyle\bullet}\quad & \{ l, m \} \text{ with } l, m \text{ in }  RR_{3t}^*
.
\end{aligned}
\end{equation*}
So for those $ d $ we have
\begin{equation} \label{Fd1}
\begin{aligned}
 F_d =  L_{3t}
& + \langle \{ f_P , f_Q \} \text{ with } 2t+1 \le \deg(Q) \le \deg(P) \le d \rangle
\\
& + 
\langle \{ f_P , l \} \text{ with } 2t+1 \le \deg(P) \le d \text{ and } l \text{ in } RR_{3t}^* \rangle
\,.
\end{aligned}
\end{equation}

We want to show in~\eqref{Fd2} that we can use simpler generators
if we work modulo~$ F_{d-1} $ instead of~$ L_{3t} $.

\begin{lemma} \label{PQlemma}
For~$ P \ne O $ the following hold.

(1)
If $ \deg(P) = d \ge 3t+1 $ and  $ 2t+1 \le \deg(Q) \le d $, then $ \{f_P, f_Q\} $ is in 
\begin{equation*}
\quad\qquad
L_{d-1} +
\langle \{f_R, l\} \text{ with } l \text{ in } RR_{d-1}^* \text{ and } 2t+1 \le \deg(R) \le d \rangle 
\,.
\end{equation*} 

(2)
If $ 2t+1 \le \deg(Q) \le \deg(P) \le 3t $, then $ \{f_P, f_Q\} $ is in
\begin{equation*}
\quad\qquad
L_{3t} +
\langle \{f_R, l\} \text{ with } l \text{ in } RR_{3t}^* \text{ and } 2t+1 \le \deg(R) \le 3t \rangle 
\,.
\end{equation*} 
\end{lemma}

\begin{proof}
(1)
Write $ d' = \deg(Q) $ and consider the linear map 
\begin{alignat*}{1}
RR_{d'-1}(D_P) \times RR_{d-1}(D_Q) \times RR_{d-1} & \to RR_{d+d'-t-1} 
\\
(l_P, l_Q, l) & \mapsto f_P l_P + f_Q l_Q - l
\,.
\end{alignat*}
The left-hand side has dimension at least $ (d'-2t) + (d-2t) + (d-t) = 2d + d' - 5t $.
Because~$ d + d' - t - 1 > 2 t - 2 $, the right-hand side has dimension $ d+d'-2t $,
hence there is a non-trivial element~$ (l_P, l_Q, l) $ in the kernel.
Then~$ l_Q \ne 0 $ because otherwise~$ l(P) = 0 $, but~$ RR_{d-1}(P) = 0 $.
If $ l_P = 0 $ then $ f_Q $ is a quotient of two elements in $ RR_{d-1}^* $ and the result is clear.
If $ l_P \ne 0 $ then the result follows from~$  0 = \{ f_P l_P , f_Q l_Q \} $ if $ l = 0 $,
and from $  0 = \{ f_P l_P / l , f_Q l_Q / l \} $ if $ l \ne 0 $.

(2) 
We again let $ d'=\deg(Q) $, but now consider the linear map
\begin{alignat*}{1}
RR_{d'-1}(D_P) \times RR_{3t-1}(D_Q) \times RR_{3t} & \to RR_{d'+2t-1} 
\\
(l_P, l_Q, l) & \mapsto f_P l_P + f_Q l_Q - l
\,.
\end{alignat*}
The right-hand side has dimension $ d'+t $ and the left-hand side has dimension at least $ (d'-2t) + t + (2t+1) = d'+t+1 $,
so that there is a non-trivial triple $ (l_P, l_Q, l) $ in the kernel.
If $ l_Q = 0 $ then $ f_P = l / l_P $ and the result is clear.
If $ l_Q \ne 0 $ then we finish the proof as for~(1).
\end{proof}

Because $ L_{d-1} \subseteq F_{d-1} $,
Lemma~\ref{PQlemma}(1) allows us to rewrite~\eqref{Fd1} for $ d \ge 3t+1 $ as
\begin{equation} \label{Fd2}
\begin{aligned}
F_d & = F_{d-1} + \langle \{ f_P , l \} \text{ with } P \ne O,\ 2t+1 \le \deg(P) \le d \text{ and } l \text{ in } RR_{d-1}^* \rangle
\\
& = F_{d-1} + \langle \{ f_P , l \} \text{ with } P \ne O, \deg(P) = d \text{ and } l \text{ in } RR_{d-1}^* \rangle
\,,
\end{aligned}
\end{equation}
and from Lemma~\ref{PQlemma}(2) we see that
\begin{equation} \label{F3t}
F_{3t} = L_{3t} + \langle \{ f_P , l \} \text{ with } P \ne O,\ 2t+1 \le \deg(P) \le 3t \text{ and } l \text{ in } RR_{3t}^* \rangle
\,.
\end{equation}
In particular, $ F_d/F_{d-1} $ for $ d \ge 3t+1 $ is generated by
the classes of symbols $ \{f_P, l\} $ for~$ P \ne O $ with $ \deg(P) = d $ and $ l $ in $ RR_{d-1}^* $.
We want to reduce their number in a given element of~$ F_d / F_{d-1} $.

\begin{lemma} \label{abdlemma}
Assume $ P \ne O $ and $ d = \deg(P) \ge t+1 $.
If $ a^\pm $ in $ RR_{e^\pm}^* $ are given, for some $ e^\pm \ge 0 $,
such that $ a^\pm(P) \ne 0 $, and there exist integers
$ \d \ge 0 $ and $ \e^\pm \le d-1 $ satisfying 
$ e^- + \e^+ \le 2d+\d-1 $, $ e^+ + \e^- \le 2d+\d-1 $
and $ \e^+ + \e^- \ge d+\d+2t-1 $,
then there exist $ b^\pm $ in $ RR_{\e^\pm}^* $ and $ g $ in $ RR_{d+t-1}(D_P) $ with
$ a^- b^+ - a^+ b^- = f_P g $.
\end{lemma}

\begin{proof}
We consider the linear map
\begin{alignat*}{1}
 RR_{\e^+} \times RR_{\e^-} \times RR_{d+t-1}(D_P) & \to RR_{2d+\d-1}
\\
 (b^+ , b^-, g) & \mapsto a^- b^+ - a^+ b^- - f_P g
\,.
\end{alignat*}
Note that the three terms are in the right-hand side because of the conditions on 
$ e^- + \e^+ $, $ e^+ + \e^- $ and $ \d $.
Because $ 2d+\d-1 > 2t-2 $, this right-hand side has dimension $ 2d+\d-t $, whereas the 
dimension of the left-hand side
is at least $ (\e^+ + 1 - t) + (\e^- + 1 - t) +(d-t) \ge 2d + \d - t + 1 $,
so that there is a non-trivial triple $ (b^+, b^-, g) $ in its kernel.
If $ b^- = 0 $ then evaluating at~$ P $ gives $ b^+(P) = 0 $, so that~$ b^+ = 0 $
because of degrees, hence $ g=0 $ as well.  Therefore $ b^- \ne 0 $, and a similar argument shows $ b^+ \ne 0 $
as well.
\end{proof}

\begin{lemma} \label{123lemma}
Suppose that $ P \ne O $, $ d = deg(P) \ge 3t+1 $, and that
$ l_1^+ $, $ l_2^+ $ and $ l_1^- $ are in $ RR_{d-1}^* $.  Then there exist $ m^\pm $ in 
$ RR_{d-1}^* $ such that 
\begin{equation*}
\{f_P,l_1^+\} + \{f_P,l_2^+\} - \{f_P,l_1^-\} \equiv  \{f_P,m^+\} - \{f_P,m^-\}
\end{equation*}
modulo $ F_{d-1} $.
\end{lemma}

\begin{proof}
We apply Lemma~\ref{abdlemma} with $ a^+ = l_1^+ $, $ a^- = l_1^- $, $ \d = 0 $, $ e^\pm = d-1 $,
$ \e^+ = 2t $ and $ \e^- = d-1 $, so we obtain  $b^+$  in  $RR_{2t}^*$, $b^-$  in  $RR_{d-1}^* $, and
$ g $ in $ RR_{d+t-1}(D_P) $ with $ l_1^- b^+ - l_1^+ b^- = f_P g $.
If $ g = 0 $ then $ 0 = \{ f_P , l_1^+ b^- /(l_1^- b^+) \} $.
If $ g \ne 0 $ then we obtain
$ 0 = \{ l_1^+ b^- / (-f_P g) , l_1^- b^+/(f_P g) \} \equiv \{ f_P , l_1^+ b^-/(l_1^- b^+) \} $ modulo
$ F_{d-1} $ because, with~$ D_P $ effective of degree~$ t \le d-1 $,
we have~$ |(g)| \subseteq C_{\le d-1}\cl $.
In either case, working modulo $ F_{d-1} $ we may replace $ l_1^+ $ with $ b^+ $ and~$ l_1^- $ with $ b^- $.
Repeating this argument with $ a^+ = l_2^+ $ and $ a^- = b^- $
we may assume that $ l_1^+ $ and $ l_2^+ $ are in $ RR_{2t}^* $, and that $ l_1^- $ is in $ RR_{d-1}^* $.

We then apply Lemma~\ref{abdlemma} with $ a^+ = l_1^+ l_2^+ $, $ a^- = l_1^- $, $ \d = t $, $ e^+ = 4t $,
$ e^- = d-1 $ and $ \e^\pm = d-1 $, obtaining $ m^\pm $ in $ RR_{d-1}^* $ and $ g $ in $ RR_{d+t-1}(D_P) $ with
$ l_1^- m^+ - l_1^+ l_2^+ m^- = f_P g $.  If $ g=0 $ then $ l_1^+ l_2^+ / l_1^- = m^+/m^- $ in $ F^* $
and the result is clear.
If $  g \ne 0 $ then
$ 0 = \{ l_1^+ l_2^+ m^-/(-f_P g) , l_1^- m^+ / (f_P g) \} \equiv \{ f_P , l_1^+ l_2^+ m^- / (l_1^- m^+) \} $
modulo $ F_{d-1} $ since again $ |(g)| \subseteq C_{\le d-1}\cl $.
\end{proof}

\begin{proof}[Proof of Theorem~\ref{introker1}]
The  maps exist because~$ T_P $ for~$ P $ in~$ C_{\ge d}\cl $ is trivial on~$ F_{d-1} $.

For the surjectivity of the first map, fix~$ P \ne O $
with~$ d = \deg(P) \ge 3 t + 1 $. Then the image
under the tame symbol of~$ \{f_P, l^+/l^- \} $
with~$ l^\pm $ in~$ RR_{d-1}^* $ is concentrated
at~$ P $ and points in~$ \{O\} \cup C_{\le d-1} $ by~\eqref{fPdef}.
Because evaluation at~$ P $ gives an injection of~$ RR_{d-1} $ to~$ k(P) $,
the surjectivity follows from Lemma~\ref{VW}.

For its injectivity, and using~$ l_P^\pm =1 $ if necessary, we
see from the last line of~\eqref{Fd2}
that~$ F_d / F_{d-1} $ is generated by the classes of~$ \{ f_P , l_P^+ \} - \{ f_P , l_P^- \} $
for~$ P \ne O $ with~$ \deg(P) = d $ and~$l^\pm $ in~$ RR_{d-1}^* $.
Using Lemma~\ref{123lemma} to collect terms for a fixed~$ P $,
we see that for any~$ P $ we need only one pair~$ l_P^\pm $ in a
given class of~$ F_d/F_{d-1} $.
The image of this class under~$ T_P $ for~$ P \ne O $ in~$ C_d $
with~$ d \ge 3 t + 1 $ is determined entirely by~$ \{ f_P , l_P^+ \} - \{ f_P , l_P^- \} $.
So, in order to have trivial image at~$ P $, we must have~$ l_P^+(P) = l_P^-(P) $ in $ k(P) $, which implies that~$ l_P^+ = l_P^- $ because of degrees.

The second isomorphism in the theorem follows easily from the
first and the fact that~$ K_2(F) = \cup_d F_d $.
The last statement follows by viewing~$ F_{d-1} $
as subgroup of~$ K_2(F) $, and~$ \coprod_{P \in \{O\} \cup C_{\le d-1}\cl} k(P)^* $
as subgroup of~$ \coprod_{P \in C\cl} k(P)^* $, and noticing
that the resulting map on the quotients is the second isomorphism
in the theorem.
\end{proof}

\section{Filtration using Riemann-Roch} \label{RRsection}

Our goal in this section is to prove Theorem~\ref{introker2}.
For this, we prove Proposition~\ref{RRprop}, which contains
a result similar to the two isomorphisms in Theorem~\ref{introker1},
but for certain subgroups of~$ F_{3t} $.
Since Theorem~\ref{introker1} already gives the classical results for~$t=0$, we now
assume that~$t \ge 1$.
Apart from that, our assumptions and notation are as in Section~\ref{degreesection}.

Let us first explain the case $t=1$, where the argument is much simpler, and does not require
the technical content of some of the lemmas below.
From Theorem~\ref{introker1} and~\eqref{F3t} we know that
the kernel of the tame symbol is contained in
\begin{equation*}
F_3 = L_3 + \langle \{ f_P , l \} \text{ with } \deg(P) = 3 \text{ and } l \text{ in } RR_3^* \rangle 
\,.
\end{equation*}
If $ \deg(P) = 3 $ and $ RR_3((P)) \ne \{0\} $ then $ f_P $ is in $ \langle RR_3^* \rangle $ by Lemma~\ref{3tlemma}(3).
(In fact, then~$ f_P $ gives a basis of~$ RR_3((P)) $.)
So we only have to consider those $ P $ in~$ C_3\cl $ such that, for~$ l $ in~$ RR_3 $,  we
have~$ l(P) = 0 $ only if~$ l = 0 $.  Let us call those $ P $ non-special.
For each such~$ P $, evaluation at~$ P $ maps~$ RR_3 $ isomorphically to~$ k(P) $,
hence~$ T_P $ maps~$ \langle \{f_P, l\} \text{ with } l \text{ in } RR_3^* \rangle
$ surjectively to~$ k(P)^* $.
Because $ T_P $ is trivial on $ L_3 $ for non-special $ P $, we concentrate on computing the
kernel of
\begin{equation} \label{L3mod}
T_P :  L_3 + \langle \{ f_P , l \} \text{ with } l \text{ in } RR_3^* \rangle / L_3 \to k(P)^* 
\end{equation}
for a fixed non-special $ P $.

Using $ l^\pm = 1 $ if necessary, the left-hand side of~\eqref{L3mod} is generated by the classes of
$ \{f_P, l^+ \} - \{f_P , l^-\} $ with $ l^\pm $ in $ RR_3^* $.
Working modulo $ L_3 $, 
we first replace such a pair with a similar pair where $ l^\pm $ in $ RR_2^* $.  For this we notice that the linear map
\begin{alignat*}{1}
RR_2 \times RR_2 \times RR_3(D_P) & \to RR_5
\\
(m^+, m^-, g) & \mapsto l^- m^+ - l^+ m^- - f_P g
\end{alignat*}
must have a non-trivial triple $ (m^+, m^-, g) $ in its kernel because of dimensions.
Evaluating $ 0 = l^- m^+ - l^- m^+ - f_P g $ at $ P $ shows~$ m^\pm(P) $
are either both zero or non-zero because~$ P $ is non-special.
But $ m^\pm(P) = 0 $ similarly implies $ m^\pm = 0 $, hence $ g=0 $ as well, which is not possible, so that $ m^\pm \ne 0 $.
When $ g = 0 $ we get $ l^+/ l^- = m^+ / m^- $ in~$ F^* $ and we can replace $ l^\pm $ with $ m^\pm $.
For $ g \ne 0 $ we have $ |(g)| \subseteq C_{\le 2}\cl $, hence~$ g $ is in $ \langle RR_3^* \rangle $
by~\eqref{fPprod} as well as parts~(1) and~(2) of Lemma~\ref{3tlemma}.
Then $  0 = \{ l^+ m^-/(-f_Pg) , l^- m^+ / (f_P g) \} \equiv \{ f_P , l^+ m^- / (l^- m^+) \} $ modulo $ L_3 $,
so working modulo $ L_3 $ we can replace $ l^\pm $ with $ m^\pm $.

We next show that for $ l_1^\pm  $ and $ l_2^+ $ in $ RR_2^* $ there exist $ m^\pm  $ in $ RR_3^* $
such that
\begin{equation} \label{t=1il}
 \{f_P , l_1^+\} +  \{f_P , l_2^+\} - \{f_P , l_1^-\} \equiv  \{f_P , m^+\} - \{f_P , m^-\}
\end{equation}
modulo $ L_3 $.  For this, with~$ f_P' $ defined after~\ref{fPdef}, we consider the linear map
\begin{alignat*}{1}
RR_3 \times RR_2 \times RR_2 & \to RR_6
\\
(m^+, m^-, g) & \mapsto l_1^- m^+ - l_1^+ l_2^+ m^- - f_P f_P' g
\,,
\end{alignat*}
notice that it has a non-trivial kernel, and, again from evaluation at $ P $,
that $ m^\pm \ne 0 $ for a non-trivial triple in its kernel.  The resulting identity gives us~\eqref{t=1il} as before
because~$ f_P' $ is in~$ RR_2^* $.
(Note that we could have allowed~$ l_1^- $ to be in $ RR_3^* $, and this type of flexibility is necessary
in the general case below.)

Since, modulo $ L_3 $, we may now replace~$ m^\pm $ again by elements in~$ RR_2^* $ by the earlier argument, repeating
the reduction in number of terms as in~\eqref{t=1il},
showing that
a class in the left-hand side of~\eqref{L3mod} is of the form~$ L_3 +  \{ f_P , l_P^+ \} - \{ f_P , l_P^- \} $
with~$ l_P^\pm $ in~$ RR_3^* $.
As noted before, the evaluation map $ RR_3 \to k(P) $ at~$ P $ is an isomorphism for $ P $ non-special.
So if the class has trivial tame symbol at the non-special~$ P $ in $ C_3\cl $,
then~$ l_P^+ = l_P^- $ and the term involving~$ f_P $ is trivial.
It then follows that~\eqref{L3mod} is an isomorphism.

In fact, these arguments work independently for all non-special~$ P $,
and show that the tame symbol~\eqref{TAME} induces an
isomorphism~$ F_3 / L_3 \rightiso \coprod_{P \in C_3\cl \setminus S_3} k(P)^* $
where~$ C_3\cl \setminus S_3$ consists of our non-special points.
(This notation agrees with~$ S $ as in Definition~\ref{Sdef} below.)
This is the first isomorphism in Proposition~\ref{RRprop} in this
case. As in the general case below, for~$ t = 1 $ one easily deduces the remainder of that proposition
from this, and, by using Theorem~\ref{introker1} for~$ d = 4 $, also Theorem~\ref{introker2}.

\medskip

We now tackle the general case, with $ t \ge 1 $.

\begin{definition} \label{Sdef}
(1)
Let $ S = \{ P \text{ in } C\cl \text{ with } RR_{3t}((P)) \ne \{0\} \} $.

(2)
For $ d=2t,\dots,3t $ we let, inside~$ K_2(F) $,
\begin{equation*}
\quad\qquad
 F_d' = L_{3t} + \langle \{f_P, l\} \text{ with } 2t+1 \le \deg(P) \le d, \ P \text{ not in } S \text{, and } l \text{ in } RR_{3t}^* \rangle
\,.
\end{equation*}
\end{definition}

We observe that $ S' \subseteq C_{\le 2t}\cl \subseteq S \subseteq C_{\le 3t}\cl $
with~$ S' $ as in Definition~\ref{Sprimedef} because~$ O $ is in~$ C_1\cl $,
with the middle inclusion a consequence of Riemann-Roch.
Note that, for~$ P $ not in $ S $ and $ l $ in $ RR_{3t} $, having $ l \ne 0 $ is equivalent with $ l(P) \ne 0 $.
We then have
\begin{equation*}
L_{3t} = F_{2t}' \subseteq F_{2t+1}' \subseteq \cdots \subseteq F_{3t}' = F_{3t}
\end{equation*}
with~$ F_{3t}' = F_{3t} $ by~\eqref{F3t} and Lemma~\ref{3tlemma}(3).

Since $ F_d'/F_{d-1}' $ for $ d = 2t+1, \dots, 3t $ is generated by the classes of elements
\begin{equation*}
\{f_P, l^+\} - \{f_P, l^- \} \text{ with } \deg(P) = d,\  P \text{ not in } S \text{, and } l^\pm \text{ in } RR_{3t}^* 
\end{equation*}
with some $ l^\pm = 1 $ if necessary, we prove some lemmas
in order to mimic the arguments above in the case $ t=1 $.

\begin{lemma} \label{flemma}
Assume that $ P $ does not lie in $ S $ and satisfies $ d = \deg(P) \ge 2t+1 $.
Let $ e = \lceil \frac{d+1}2 \rceil + t -1 $, and assume $ \e $ satisfies $ e+1 \le \e \le 3t $.
If $ a^\pm $ are in $ RR_\e^* $, then there exist $ b^\pm $ in $ RR_{\e-1}^* $
and $ g $ in $ RR_{3t}(D_P) $, such that
$ a^- b^+ - a^+ b^- - f_P g = 0 $.
In particular,
\begin{equation*}
\{ f_P, a^+ \} - \{ f_P, a^- \}
\equiv
\{ f_P, b^+ \} - \{ f_P, b^- \}
\end{equation*}
modulo $L_{3t}  $.
\end{lemma}

\begin{proof}
We consider the linear map
\begin{alignat*}{1}  
RR_{\e-1} \times RR_{\e-1} \times RR_{3t}(D_P) & \to RR_{2\e-1}
\\
 (b^+ , b^-, g) & \mapsto a^- b^+ - a^+ b^- - f_P g
\,.
\end{alignat*}
The term $ f_P g $ is in the right-hand side because $ 3t+d-t \le 2\e-1 $, i.e., $ e+1 \le \e $.
Here~$ RR_{2\e-1} $ has dimension $ 2\e-t $ because $ \e \ge e+1 \ge t $, and the left-hand side has dimension
at least $ 2(\e-t) + (t+1) $, so there is a non-zero triple $ ( b^+, b^-, g) $ in the kernel.
With~$ P $ not in~$ S $ we have~$ a^\pm(P) \ne 0 $, and $ b^\pm = 0 $  only if $ b^\pm(P) = 0 $.
Evaluating the identity $ a^- b^+ - a^+ b^- = f_P g $ at $ P $, one
finds that~$ b^\pm \ne 0 $.
If~$ g=0 $ then $ a^+ / a^- = b^+/b^- $ in $ F^* $ and the result is clear,
whereas for~$  g \ne 0 $ we have $ 0 = \{ a^+ b^-/(-f_P g) , a^- b^+ / (f_P g) \} \equiv \{ f_P , a^+ b^- / (a^- b^+) \} $
modulo $ L_{3t} $.
\end{proof}

\begin{lemma} \label{llmmlemma}
Suppose that~$ P $ is not in~$ S $, $ d = \deg(P)  $ satisfies $ 2t+1 \le d \le 3t $
and~$ l^\pm $ are in~$ RR_{d-1}^* $. Then there exist
$ m^+ $ in~$ RR_{2t-1}^* $, $ m^- $ in~$ RR_{3t}^* $, and~$ g $
in~$ RR_{3t} $ with~$ l^- m^+ - l^+ m^- - f_P g = 0 $. In particular,
\begin{equation*}
 \{f_P , l^+\} - \{f_P , l^-\} \equiv  \{f_P , m^+\} - \{f_P , m^-\}
\end{equation*}
modulo $ F_{d-1}' $.
\end{lemma}

\begin{proof}
We consider the linear map
\begin{alignat*}{1}  
RR_{2t-1} \times RR_{3t} \times RR_{d+t-1}(D_P) & \to RR_{d+3t-1}
\\
 (m^+ , m^-, g) & \mapsto l^- m^+ - l^+ m^- - f_P g
\,.
\end{alignat*}
Note that $f_P g$ is in the right-hand side because $ 2d-1 \le d+3t-1 $.
The dimension on the right-hand side is $d+2t$, whereas the dimension on the left-hand side is at least
$ t + (2t+1) + (d-t) = d+2t+1 $ so that the kernel is non-trivial.
One checks as before that, for a non-trivial triple in the kernel, $m^\pm\ne0$ because $P$ is not in $ S$.
If~$ g=0 $ then~$ l^+/l^- = m^+/m^- $ in $ F^* $ and the result
is clear.
If~$ g \ne 0 $ then~$ |(g)| \subseteq C_{\le d-1}\cl $ because~$ D_P \ge 0 $
and has degree~$ t $.
It follows from~\eqref{fPprod}, together with all parts of Lemma~\ref{3tlemma},
that~$ g $ is a product of an element in~$ \langle RR_{3t}^* \rangle $
and functions~$ f_R $ for~$ R $ not in~$ S $ and~$ \deg(R) = 2t+1, \dots, d-1 $.
It follows that in this case
$ 0 = \{ l^+ m^-/(-f_Pg) , l^- m^+/(f_Pg) \} \equiv \{ f_P, l^+ m^- /(l^- m^+) \}$ modulo~$ F_{d-1}' $.
\end{proof}

\begin{lemma} \label{3t3t3tlemma}
Suppose that $ P $ is not in $ S $, that $ d = \deg(P) $ satisfies $ 2t+1 \le d \le 3t $, 
that $ a_1^+ $ and~$ a_2^+ $ are in~$ RR_{2t-1}^* $,
and that $ a_1^- $ is in $ RR_{3t}^* $.  Then there exist~$ b^\pm $ in~$ RR_{3t}^* $ and~$ g $ in~$ RR_{3t} $
satisfying $ a_1^- b^+ - a_1^+ a_2^+ b^- - f_P f_P' g = 0 $.
In particular,
\begin{equation*}
\{ f_P, a_1^+ \} + \{ f_P, a_2^+ \} - \{ f_P, a_1^- \}
\equiv
\{ f_P, b^+ \} - \{ f_P, b^- \}
\end{equation*}
modulo $ L_{3t} $.
\end{lemma}

\begin{proof}
We consider the linear map
\begin{alignat*}{1}  
RR_{3t} \times RR_{3t} \times RR_{3t} & \to RR_{7t}
\\
 (b^+ , b^-, g) & \mapsto a_1^- b^+ - a_1^+ a_2^+ b^- - f_P f_P' g
\,.
\end{alignat*}
As we saw after~\eqref{fPdef}, we have~$ f_P f_P' $ in $ RR_{d+t}((P)+D_P')^* \subseteq RR_{4t}^* $ here, so all terms are in
the right-hand side, which has dimension $ 6 t + 1 $.  The left-hand side has dimension at least $ 3(2t+1) $
so the kernel contains a non-trivial triple~$ (b^+, b^-, g) $.
As before, evaluating at~$ P $ implies that~$ b^\pm \ne 0 $ because~$ P $ is not in~$ S $.
With~$ f_P' $ in~$ RR_{2t}^* $,
the triple implies the result in the now standard way.
\end{proof}

\begin{lemma} \label{N3t3t3t}
Suppose that $ P $ is not in $ S $, that $ d = \deg(P) $ satisfies $ 2t+1 \le d \le 3t $, 
and that $ l_1^\pm $ and $ l_2^+ $ are in $ RR_{3t}^* $.
Then there exist $ m^\pm $ in $ RR_{3t}^* $ such that
\begin{equation*}
\{ f_P, l_1^+ \} + \{ f_P, l_2^+ \} - \{ f_P, l_1^- \}
\equiv
\{ f_P, m^+ \} - \{ f_P, m^- \}
\end{equation*}
modulo $ F_{d-1}' $.
\end{lemma}

\begin{proof}
Iterating Lemma~\ref{flemma} tells us that, modulo $ L_{3t} \subseteq F_{d-1}'$,
we may assume $ l_1^+ $ and $ l_1^- $ are in $ RR_e^* $ with
$ e = \lceil \frac{d+1}2 \rceil +t -1 $.  Since $ e \le d-1 $ for $ d \ge 2t+1 $, 
Lemma~\ref{llmmlemma} then tells us that, working modulo $ F_{d-1}' $,
we may assume $ l_1^+ $ is in $ RR_{2t-1}^* $ and $ l_1^- $ is in $ RR_{3t}^* $.
Applying the same procedure to $ l_2^+ $ and the new $ l_1^- $ we see that we may assume that $ l_1^+ $
and $ l_2^+ $ are in $ RR_{2t-1}^* $ and that $ l_1^- $ is in $ RR_{3t}^* $.
We then apply Lemma~\ref{3t3t3tlemma} with $ a_1^+ = l_1^+$, $ a_2^+ = l_2^+ $, and $ a_1^- = l_1^- $, 
obtaining~$ m^\pm = b^\pm $ in $ RR_{3t}^* $ for which the statement
of the lemma holds.
\end{proof}

\begin{proposition} \label{RRprop}
For $ d = 2t+1,\dots,3t $, the tame symbol induces isomorphisms
\begin{equation*}
F_d'/F_{d-1}' \rightiso \coprod_{P \in C_d\cl \setminus S_d} k(P)^*
\text{ and }
\quad
K_2(F) / F_{d-1}' \rightiso \coprod_{P \in C_{\ge d}\cl \setminus S_{\ge d}} k(P)^*
\,.
\end{equation*}
\end{proposition}

\begin{proof}
The maps exist because the image of an element of~$ F_{d-1}' $
under the tame symbol is supported in~$ S \cup C_{\le d-1}\cl $.

For the surjectivity of the first map, fix~$ P $ in~$ C_d\cl \setminus S_d $.
Then evaluation at~$ P $ gives an injection of~$ RR_{3t} $ into~$ k(P) $.
Using Lemma~\ref{VW} one sees that there are $ l $ and~$ m $
in~$ RR_{3t}^* $ such that $ \{f_P, l/m\} $ has prescribed image
in~$ k(P)^* $, and trivial image at any other point of~$ C_d\cl \setminus S_d $.
For its injectivity we note that, by Lemma~\ref{N3t3t3t}, 
any element of~$ F_d' / F_{d-1}' $ can be written as a sum
of classes~$ F_{d-1}' + \{f_Q,l_Q^+\} - \{f_Q,l_Q^-\} $
with $ l_Q^\pm $ in $ RR_{3t}^* $ for different~$ Q $ in~$ C_d \setminus S_d $.
Then~$ T_P $ being trivial on this sum for all~$ P $ in~$ C_d \setminus S_d $
implies, by taking~$ P = Q $, that~$ l_Q^+ = l_Q^- $ for each~$ Q $.

The second isomorphism follows by using~$ F_{3t}' = F_{3t} $,
the second isomorphism in Theorem~\ref{introker1} for~$ d = 3 t + 1 $
(noting that~$ S \subseteq C_{\le 3t}\cl $),
and repeated application of the first isomorphism
of the proposition for descending~$ d $.
\end{proof}

\begin{proof}[Proof of Theorem~\ref{introker2}]
The first statement is the second statement of Proposition~\ref{RRprop}
for~$ d = 2 t + 1 $ because~$ F_{2t}' = L_{3t} $ and~$ C_{\le 2t}\cl \subseteq S $.
The second statement follows from this
as in the proof of Theorem~\ref{introker1} in Section~\ref{degreesection}.
\end{proof}

\begin{remark} \label{sharpremark}
If $ t=1 $ then Theorem~\ref{introker2}
is sharp in the following sense.
If~$ E $ is an elliptic curve with $ k $-rational point~$ O $ as neutral
element for the group law, then the elements of~$ L_2 $ are invariant
under multiplication by~$ -1 $ on~$ E $, whereas there are many
examples of such~$ E $ in which the kernel of~\eqref{TAME} contains
anti-invariant elements of infinite order. Therefore, one cannot
replace~$ L_3 $ with~$ L_2 $.
If~$ k = \Q $ then the curve is modular, hence this always applies
by a theorem of Beilinson (see, e.g., \cite[Th\'eor\`eme~9.1]{scha92}).

Note also that~$ F_d = L_d $ for each~$ d $ if~$ E(k) = \{O\} $: in~\eqref{fPprod} we have~$ f' $
in~$k^* $ as~$ E_{\le1}\cl = \{O\} $, 
and $ D_P = (O) $ for every $ P \ne O $ with $ d= \deg(P) \ge 2 $, so $ f_P $ is in~$ RR_d^* $.
For such curves~$ S = C_{\le 3}\cl $, so Proposition~\ref{RRprop} and Theorem~\ref{introker2}
do not improve on Theorem~\ref{introker1}.
\end{remark}

\bibliographystyle{plain}

\end{document}